# Design, Modeling, and Geometric Control on SE(3) of a Fully-Actuated Hexarotor for Aerial Interaction

Ramy Rashad, Petra Kuipers, Johan Engelen and Stefano Stramigioli

*Abstract*— In this work we present the optimization-based design and control of a fully-actuated omnidirectional hexarotor. The tilt angles of the propellers are designed by maximizing the control wrench applied by the propellers. This maximizes (a) the agility of the UAV, (b) the maximum payload the UAV can hover with at any orientation, and (c) the interaction wrench that the UAV can apply to the environment in physical contact. It is shown that only axial tilting of the propellers with respect to the UAV's body yields optimal results. Unlike the conventional hexarotor, the proposed hexarotor can generate at least 1.9 times the maximum thrust of one rotor in any direction, in addition to the higher control torque around the vehicle's upward axis. A geometric controller on SE(3) is proposed for the trajectory tracking problem for the class of fully actuated UAVs. The proposed controller avoids singularities and complexities that arise when using local parametrizations, in addition to being invariant to a change of inertial coordinate frame. The performance of the controller is validated in simulation.

## I. INTRODUCTION

Unmanned aerial vehicles (UAVs) have become a successful cost-effective tool for various civilian applications such as surveillance, visual inspection, precision agriculture, and mapping. Many of these remote-sensing applications are already operational and economically successful. Recently there has been a growing interest in the robotics research community in using aerial vehicles for physical interaction. Thus, extending the range of applications which UAVs could be used for, such as human interaction, object manipulation or non-destructive inspection tasks.

Most multi-rotor UAVs currently used, such as quadrotors or hexarotors, are designed and optimized for flying without physical interaction with the environment. Therefore, these traditional UAVs have all their rotors aligned in a single plane with the thrust direction pointing upwards. This results in an underactuated vehicle with coupling between the translational and the rotational dynamics of the UAV, which limits the UAV's ability to track arbitrary position and attitude trajectories. In addition, this also limits the ability to perform complex physical interaction tasks with the environment as it cannot reject arbitrary reaction wrenches from the environment.

This work has been funded by the cooperation program INTERREG Deutschland-Nederland as part of the SPECTORS project number 143081.

R. Rashad, J. Engelen and S. Stramigioli are with the Robotics and Mechatronics group, Faculty of Electrical Engineering, Mathematics and Computer Science, University of Twente, Enschede, The Netherlands. Email: {r.a.m.rashadhashem,j.b.c.engelen, s.stramigioli}@utwente.nl

P. Kuipers is with the Precision Engineering chair, Faculty of Engineering Technology. University of Twente, Enschede, The Netherlands. Email: pkuipers11@gmail.com

Many research groups have studied several approaches to overcome the problem of underactuation of traditional multi-rotor UAVs for physical interaction tasks. One approach is to equip the UAV with a multi-degree of freedom (DoF) manipulator as in [1], [2]. However, the UAV/manipulator system suffers from limited payload capabilities and operation time. Furthermore, it is unable to exert arbitrary wrenches in an arbitrary direction, since the UAV itself cannot hover holding arbitrary orientations.

This has led researchers to investigate novel multi-rotor UAV designs to develop fully actuated vehicles. One of the early attempts was by adding four lateral rotors below each vertical one [3]. In addition to the decrease of payload capabilities, this design suffers from high aerodynamic interaction between the vertical and lateral rotors which complicates the flight controller design. In [4], an overactuated quadrotor was designed by tilting each propeller with a servo motor independently. The modified quadrotor was shown capable of hovering at nonzero pitch and roll angles [5], however the maximum achievable rotation angle in hover was limited.

Recently, optimization-based UAV designs have been proposed in the literature. In [6], a modified hexarotor was designed by tilting all propellers by a fixed angle about one axis. In [7]–[9], this approach was generalized by allowing the propellers of a hexarotor to have arbitrary orientations parametrized by two fixed angles. In [7], these fixed angles were optimized to increase dynamic manipulability and maximum translational acceleration, while in [8], [9], the angles were optimized for minimizing the control effort along a specific desired trajectory. In [10], an omnidirectional UAV has been designed which comprises eight rotors attached to the vertices of a cube with different orientations. The design was optimized to maximize the propellers' generated wrench while achieving dynamical properties which are almost independent from its orientation. Similarly in [11], another omnidirectional UAV was designed for maximizing generated wrench but with six rotors. The designs of both [10] and [11] are capable of independently controlling the force and torque generated due to the use of reversible propellers, which can generate both positive and negative thrust.

In this paper, we introduce the design, modeling and control of a fully-actuated omnidirectional hexarotor. The proposed hexarotor has propellers pointing in dissimilar directions similar to the design of [8]. However, the optimization criteria used in this paper is maximizing the generated wrench of the UAV and the propellers are reversible. Furthermore, in comparison to the designs of [10],

[11], the proposed design is mechanically simpler because it only comprises of a minor mechanical modification to the traditional hexarotor UAV.

In addition to the optimization-based design, we also approach the modeling and control design of the proposed UAV in a geometric approach by exploiting the configuration space $SE(3)$ of the UAV. In this paper, a generic geometric tracking controller is proposed for the class of fully-actuated UAVs. In previous research (in the aerial robotics literature) that applies geometric controllers, [11]–[13], the rotational and translational controllers are designed separately. These previous works yield controllers that are not invariant to changes in the inertial coordinate frame, in contrast to the proposed controller. Moreover, the proposed controller does not suffer from singularities or complexities that arise when using local coordinates. The approach used in this paper is based on the framework introduced by Bullo and Murray [14] for geometric tracking on $SE(3)$. Although the control objective in this work is trajectory tracking without interaction, the proposed controller is a step towards developing geometric coordinate-free interaction control strategies.

The paper is organized as follows: section II will present the optimization-based design of the proposed UAV. In section III, the dynamic model of the class of fully actuated aerial vehicles will be presented. The geometric tracking controller of this class of systems, designed on the group of rigid motions $SE(3)$, will be presented in section IV. In section V, illustrative simulation results showing the performance of the proposed controller will be presented. Finally, the paper is concluded in section VI.

## II. UAV Design

### A. Vehicle Description

A conventional hexarotor consists of six parallel propellers placed at the vertices of a planar hexagon. Similar to the quadrotor, the hexarotor can only generate thrust normal to the rotating plane of the propellers. To modify a traditional hexarotor in order to be fully-actuated, the six rotors should not be pointing in the same direction. In this work, each rotor's orientation is fixed and parametrized by two fixed angles ($\alpha_i$ and $\beta_i$), described in the next section. The tilting angles ($\alpha_i$, $\beta_i$) are set to maximize the static generated wrench of the propellers on the UAV's body. Consequently, this maximizes the UAV's agility, payload, and interaction wrench with the environment. Another requirement for an omnidirectional vehicle with fixedly-oriented propellers, is the ability of the rotors to produce thrust bi-directionally. In practice, this can be achieved by variable pitch propellers as in [15] or reversible electronic speed controllers as in [10].

Active tilting of the rotors, as in [5], is not considered in this work because it requires additional actuators for tilting, thus reducing the payload of the UAV, reducing the interaction wrench it can apply to the environment, and increasing the control system complexity. Although a reallocation of the rotors, as in [10], [11], provides more design flexibility and perhaps a higher generated wrench in the optimization process, in our design the rotors are placed at the vertices of a planar hexagon, similar to a conventional hexarotor. The reason behind this choice is that a reallocation of the rotors also increases the design complexity. This causes the designed UAV to lose one of the main features that made conventional multi-rotor vehicles popular, namely their mechanical simplicity. Moreover, the proposed modification is easily applicable to conventional hexarotors.

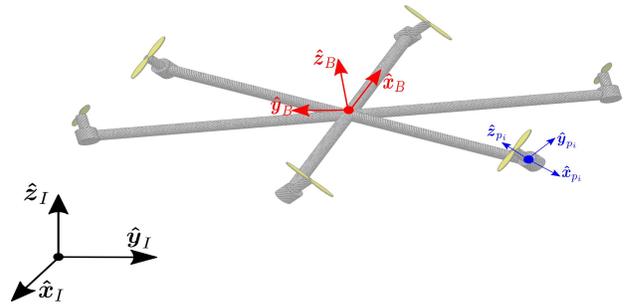

Fig. 1: Schematic view of the reference frames used.

### B. Static Wrench Analysis

Let $\{\Psi_I : o_I, \hat{\boldsymbol{x}}_I, \hat{\boldsymbol{y}}_I, \hat{\boldsymbol{z}}_I\}$ denote a right handed orthonormal inertial frame and $\{\Psi_B : o_B, \hat{\boldsymbol{x}}_B, \hat{\boldsymbol{y}}_B, \hat{\boldsymbol{z}}_B\}$ denote a body fixed frame attached to the UAV's center of mass (CoM) and aligned with the principal inertia axes of the UAV, as shown in Fig. 1. Let $\{\Psi_{p_i} : o_{p_i}, \hat{\boldsymbol{x}}_{p_i}, \hat{\boldsymbol{y}}_{p_i}, \hat{\boldsymbol{z}}_{p_i}\}$ denote the frame associated with the $i$-th propeller, where $\hat{\boldsymbol{z}}_{p_i}$ is the direction of generated thrust and the origin $o_{p_i}$ coincides with the CoM of the $i$-th propeller. The attaching location of the $i$-th propeller, for $i \in \mathcal{N}_p := \{1, \cdots, 6\}$ in $\Psi_B$ is given by

$$\boldsymbol{r}_i := \boldsymbol{o}_{p_i}^B = \boldsymbol{R}_z(\psi_i)[L, 0, 0]^\top, \qquad (1)$$

where $\boldsymbol{R}_z(\cdot) \in SO(3)$ is a rotation matrix about the $z$ axis, $L$ is the distance from the hexarotor's central axis $\hat{\boldsymbol{z}}_B$ to each rotor, and the angle $\psi_i := (i-1)\frac{\pi}{3}$. The orientation of $\Psi_{p_i}$ with respect to $\Psi_B$ is given by

$$\boldsymbol{R}_{p_i}^B = \boldsymbol{R}_z(\psi_i)\boldsymbol{R}_x(\alpha_i)\boldsymbol{R}_y(\beta_i), \quad i \in \mathcal{N}_p, \qquad (2)$$

where the angles $\alpha_i$ and $\beta_i$ uniquely define the direction of the thrust generation axis $\hat{\boldsymbol{z}}_{p_i}$ in $\Psi_B$ as shown in Fig. 2.

It is well known from literature that the aerodynamic thrust and drag torque of a single propeller are both proportional to the square of the propeller's spinning velocity [13]. The thrust magnitude generated by the $i$-th propeller in $\Psi_{p_i}$ will be denoted by $\lambda_i \in \mathbb{R}$, while the drag torque $\tau_{d,i}$ will be expressed as

$$\tau_{d,i} = \gamma \sigma_i \lambda_i, \qquad (3)$$

where $\gamma$ is the propeller-specific drag-to-thrust ratio, and $\sigma_i \in \{-1, 1\}$ specifies the direction of the propeller's rotation.

From the aforementioned definitions, the cumulative propellers' generated wrench in $\Psi_B$ can be written as

$$\boldsymbol{W}_p^B = \begin{bmatrix} \boldsymbol{\tau}_p^B \\ \boldsymbol{f}_p^B \end{bmatrix} = \begin{bmatrix} \sum_i \lambda_i(\boldsymbol{r}_i \wedge \boldsymbol{u}_i + \gamma \sigma_i \boldsymbol{u}_i) \\ \sum_i \lambda_i \boldsymbol{u}_i \end{bmatrix}, \qquad (4)$$

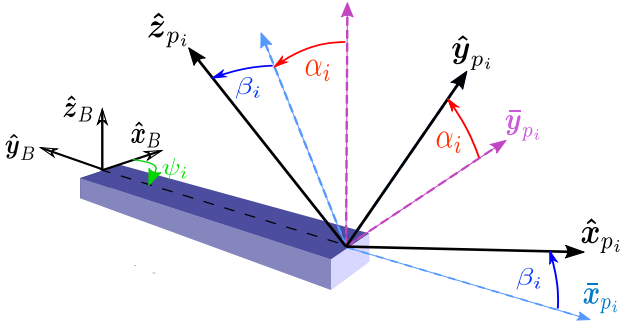

Fig. 2: The orientation of $\Psi_{p_i}$ with respect to $\Psi_B$ parametrized by $\psi_i, \alpha_i$ and $\beta_i$. First a rotation of $\psi_i$ about $\hat{z}_B$, then a rotation of $\alpha_i$ about $\bar{x}_{p_i}$, and then a rotation $\beta_i$ about $\bar{y}_{p_i}$.

where $\wedge$ denotes the vector product in $\mathbb{R}^3$, $\boldsymbol{u}_i = \boldsymbol{R}^B_{p_i}\hat{\boldsymbol{e}}_3$ denotes the thrust generation direction of the $i$-th propeller, $\hat{\boldsymbol{e}}_i$ is a vector of zeros with one at the $i$-th element, and $\boldsymbol{r}_i$ is given in (1) (with norm $\|\boldsymbol{r}_i\| = L$).

Because the drag torque is typically an order of magnitude smaller than the remaining terms [10], it will be neglected during the optimization process. Another reason for neglecting the drag torque is to allow the scalability of the optimization results ( validated in the subsequent section). Therefore, the static generated wrench of the propellers on the UAV's body can be expressed as

$$\boldsymbol{W}^B_p = \begin{bmatrix} \boldsymbol{\tau}^B_p \\ \boldsymbol{f}^B_p \end{bmatrix} = \begin{bmatrix} \boldsymbol{t}_1, \cdots, \boldsymbol{t}_6 \\ \boldsymbol{u}_1, \cdots, \boldsymbol{u}_6 \end{bmatrix} \boldsymbol{\lambda} = \begin{bmatrix} \boldsymbol{T} \\ \boldsymbol{U} \end{bmatrix} \boldsymbol{\lambda} = \boldsymbol{M}\boldsymbol{\lambda}, \quad (5)$$

where $\boldsymbol{\lambda} = [\lambda_1, \cdots, \lambda_6]^\top$, $\boldsymbol{t}_i = \boldsymbol{r}_i \wedge \boldsymbol{u}_i + \gamma \sigma_i \boldsymbol{u}_i$, and $\gamma$ will be set to be zero during the optimization process presented next.

### C. Optimization-based design

In this section, we discuss the optimization of the fixed tilting angles of the propellers $\alpha_i$ and $\beta_i$. The objective of the optimization process is to maximize the aerodynamic wrench that can be generated by the propellers on the body of the UAV. Another design requirement is having a fully actuated UAV which can independently control its force and torque in arbitrary directions. This is achieved when the mapping $\boldsymbol{M}$ (Eq. 5) between the propellers' thrust vector $\boldsymbol{\lambda}$ and the generated wrench $\boldsymbol{W}^B_p$ has full rank.

Based on the analysis in [7], [8], it can be concluded that a constrained design of the twelve design parameters ($\alpha_i$ and $\beta_i$, for $i = 1, \cdots, 6$) is advantageous. This yields no asymmetries in the UAV body which helps reduce the coupling between the force and torque generated [11]. In addition, this will cause a minimum change to the vehicle's CoM from the UAV's geometric center. For the proposed UAV design, we constrain $\alpha_i$ and $\beta_i$ to be

$$\alpha_1 = -\alpha_2 = \alpha_3 = -\alpha_4 = \alpha_5 = -\alpha_6 = \alpha, \quad (6)$$
$$\beta_1 = -\beta_2 = \beta_3 = -\beta_4 = \beta_5 = -\beta_6 = \beta. \quad (7)$$

It is shown in [8] that the aforementioned choice yields minimum control effort for a given trajectory. Thus, the design parameters are now decreased to two, $\alpha$ and $\beta$, which will be designed next for maximizing the generated wrench.

To formulate the optimization problem, we introduce the following sets: let $\mathcal{CF}$ and $\mathcal{CT}$ denote the set of control forces and control torques, respectively, defined by

$$\mathcal{CF} := \{\boldsymbol{f} \in \mathbb{R}^3 | \boldsymbol{f} = \boldsymbol{U}\boldsymbol{\lambda}, \boldsymbol{\lambda} \in \Lambda\}, \quad (8)$$
$$\mathcal{CT} := \{\boldsymbol{\tau} \in \mathbb{R}^3 | \boldsymbol{\tau} = \boldsymbol{T}\boldsymbol{\lambda}, \boldsymbol{\lambda} \in \Lambda\}, \quad (9)$$

where $\Lambda$ is the set of feasible propellers' thrusts given by

$$\Lambda := \{\boldsymbol{\lambda} \in \mathbb{R}^6 | \quad \|\boldsymbol{\lambda}\|_\infty \leq \lambda_{\max}\}, \quad (10)$$

where $\lambda_{max}$ is the maximum thrust output of each rotor, and $\|\cdot\|_\infty$ denotes the infinity-norm. As discussed earlier, it is assumed that the rotors produce bi-directional thrust and that all rotors have equivalent propellers and thrust capabilities.

The optimized control force and torque sets defined above are visualized in Fig. 3. These sets encompass all the feasible forces and torques that the UAV can generate given its rotors' locations and orientations, encoded in the mappings $\boldsymbol{U}$ and $\boldsymbol{T}$. The control force and control torque sets of a multi-rotor configuration can be formed by examining the positive range space of the mappings $\boldsymbol{U}$ and $\boldsymbol{T}$ respectively [16]. A vehicle with high wrench generation capability would have a larger volume of control forces and torques than a vehicle with low generated wrench capability. Thus, by changing the rotor locations and/or orientations, these sets can be maximized. In this study, the rotor locations are predetermined and only the rotor orientations, parametrized by $\alpha$ and $\beta$, need to be designed.

To mathematically represent the size of the control force and torque sets, we define the minimum guaranteed control force and control torque for all orientations to be

$$\mathcal{F}_{\min} := \min_{i,j \in \mathcal{N}_p} \lambda_{\max} \sum_{k \in \mathcal{N}_p} \frac{|(\boldsymbol{u}_i \wedge \boldsymbol{u}_j)^\top \boldsymbol{u}_k|}{\|\boldsymbol{u}_i \wedge \boldsymbol{u}_j\|}, \quad (11)$$

$$\mathcal{T}_{\min} := \min_{i,j \in \mathcal{N}_p} \lambda_{\max} \sum_{k \in \mathcal{N}_p} \frac{|(\boldsymbol{t}_i \wedge \boldsymbol{t}_j)^\top \boldsymbol{t}_k|}{\|\boldsymbol{t}_i \wedge \boldsymbol{t}_j\|}. \quad (12)$$

The derivation of the aforementioned expressions can be found in [11] based on a similar analysis done for cable-driven robots in [16].

The parameters $\mathcal{F}_{\min}$ and $\mathcal{T}_{\min}$ are interpreted geometrically as the radii of the maximum sphere that can be contained within $\mathcal{CT}$ and $\mathcal{CF}$ respectively, as illustrated in Fig. 3. Recall the definitions of $\boldsymbol{u}_i$ and $\boldsymbol{t}_i$ below equations (4) and (5) respectively; the minimum guaranteed force $\mathcal{F}_{\min}$ depends on the rotors' orientations only, while the minimum guaranteed torque $\mathcal{T}_{\min}$ depends on both the orientations and the locations of the rotors. In the optimization-based design by [11], the rotors' orientations was designed first by maximizing $\mathcal{F}_{\min}$ then the rotors' locations were determined by maximizing $\mathcal{T}_{\min}$. However, in this study, the rotors' locations are predetermined. Thus, the angles $\alpha$ and $\beta$ parametrizing the rotors' orientations can be designed by maximizing $\mathcal{F}_{\min}$, $\mathcal{T}_{\min}$, or a weighed combination.

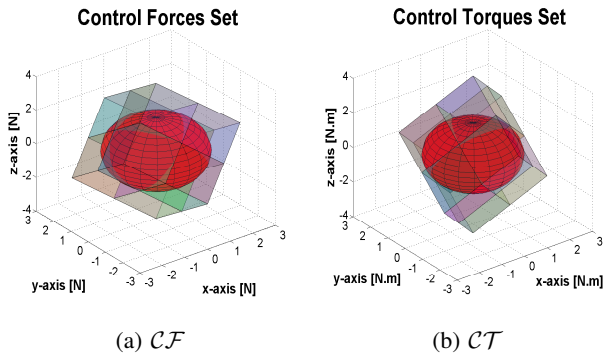
(a) $\mathcal{CF}$

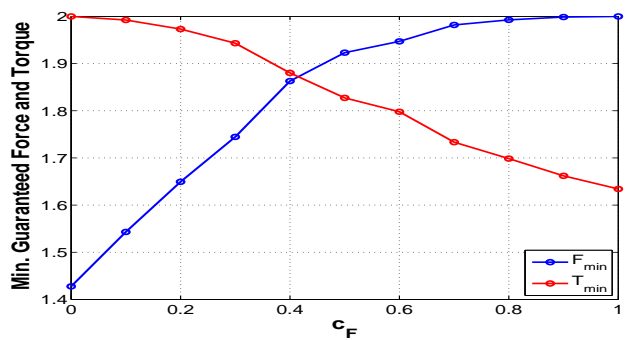
(b) $\mathcal{CT}$

Fig. 3: Control forces and torques sets. Computed for optimal $\alpha$ and $\beta$ with $L=1, \gamma=0, \lambda_{max}=1$. Within each set is a sphere of radius $\mathcal{F}_{\min}$ and $\mathcal{T}_{\min}$, respectively.

Now we can formulate the design problem as a constrained optimization problem defined as

$$\max_{\alpha,\beta} \; \zeta(\alpha,\beta) := c_F \mathcal{F}_{\min}(\alpha,\beta) + \frac{(1-c_F)}{L}\mathcal{T}_{\min}(\alpha,\beta) \quad (13)$$
$$\text{subj. to} \quad 0 \le \alpha \le \frac{\pi}{2}, \quad 0 \le \beta \le \frac{\pi}{2},$$

where $0 \le c_F \le 1$ is the weighting factor. In order for the design process to be generic and independent of the size of the UAV or the maximum thrust of the rotors, it will be assumed that the actuators lie on a unit sphere, i.e. $L=1$, and that the maximum thrust $\lambda_{max}=1$. Moreover, due to the normalization of $\mathcal{T}_{\min}$ in (13) by $L$, even if dimensionally consistent, it will be shown that any other choice of $L$ will not change the optimization results (for $\gamma=0$).

### D. Optimization Results

In this section the results of the optimization problem stated in (13) will be presented and discussed. We will discuss three different cases of the optimization-process corresponding to three choices of the weighting factors. In Case 1, the design parameters are optimized by maximizing $\mathcal{F}_{\min}$ only ($c_F=1$). In Case 2, the design parameters are optimized by maximizing $\mathcal{T}_{\min}$ only ($c_F=0$). While in Case 3, the design parameters are optimized by maximizing both $\mathcal{F}_{\min}$ and $\mathcal{T}_{\min}$ with equal weights ($c_F=0.5$). The general effect of the weighing factor on the optimized values of $\mathcal{F}_{\min}$ and $\mathcal{T}_{\min}$ is displayed in Fig. 4.

Due to the fact that the search space is parametrized by two variables, the cost function can be visualized as a two-dimensional surface, as shown in Fig. 5. For Case 1, it can be observed that there is no unique optimum for the maximum value $\mathcal{F}_{\min}=2$. Thus, several combinations of the optimal angles $\alpha^*$ and $\beta^*$ could be chosen. These combinations also include the two extremes, $\alpha^*=0$ for some $\beta^* \ne 0$, and $\beta^*=0$ for some $\alpha^* \ne 0$. Thus, it can be concluded that for the first choice of the cost function, designing two tilting angles is redundant and only one of the two is sufficient. For Case 2, it can be observed that the maximum $\mathcal{T}_{\min}$ is achieved only for $\beta^*=0$. The same can be concluded for Case 3, as well as any other case with $c_F < 1$.

Fig. 4: Effect of the weighing factor $c_F$ on the optimized minimum guaranteed force and torque.

TABLE I: Minimum guaranteed control force and torque values for different optima.

| Optimum w.r.t | $\alpha^*$ (°) | $\mathcal{F}_{\min}$ | $\mathcal{T}_{\min}$ |
|---|---|---|---|
| Case 1: $c_F=1$ | 54.7 | 2 | 1.633 |
| Case 2: $c_F=0$ | 35.3 | 1.414 | 2 |
| Case 3: $c_F=0.5$ | 47.7 | 1.912 | 1.838 |

Therefore, it is concluded that having a tilting angle $\beta=0$ in all cases is advantageous, not to mention the additional mechanical simplicity in physically realizing the design. For the computed spaces and $\beta=0$, the optimal $\alpha^*$ and the corresponding cost values are presented in Table I.

Table II lists the differences between the three optimization cases. The highest $\mathcal{F}_{\min}$ is obtained in the first case, likewise the highest $\mathcal{T}_{\min}$ is obtained in the second case. Both $\mathcal{F}_{\min}$ and $\mathcal{T}_{\min}$ decrease for the combined criterion because a trade off has been made. However, as can be seen from the table, for a 4% decrease in $\mathcal{F}_{\min}$, $\mathcal{T}_{\min}$ is increased by 12% with respect to the first case. With respect to the second case, the $\mathcal{T}_{\min}$ is only decreased by 8%, but $\mathcal{F}_{\min}$ is increased by 35%. Thus, it can be concluded that the $\alpha^*$ obtained for the combined criterion will be regarded as the optimal solution, used in the visualization of $\mathcal{CF}$ and $\mathcal{CT}$ in Fig. 3.

In conclusion, for interaction purposes, it is sufficient to design the hexarotor with only axial tilting ($\alpha$). Compared to the conventional hexarotor, the proposed UAV has a significant advantage: while the conventional hexarotor can generate a force 6 times the maximum thrust $\lambda_{max}$ in the upward direction ($\hat{z}_B$) only, the proposed hexarotor can generate forces at least 1.9 times the maximum thrust in all directions. The force generated in the upward direction can reach up to 4 times the maximum thrust, as seen in Fig. 3. Regarding the control torque about the upward axis, the conventional hexarotor is limited to only the torque from the drag torque, whereas the proposed UAV has significantly higher control torque available.

*Remark* 1: In the optimization process presented, it was assumed that $L=1$ m and $\gamma=0$. As shown in Fig. 6, the optimum $\alpha^*$ is independent of the choice of $L$, due to the normalization of $\mathcal{T}_{\min}$ in (13) by $L$. However, it can be observed in Fig. 6 that for $\gamma \ne 0$, the optimum is asymptotically constant for high values of $L$. Whereas, for $L <= 0.3$ m the optimum depends dramatically on the value of $\gamma$. Thus, for small scale UAVs, it is recommended that the optimization process be performed for a specific choice of $\gamma$.

TABLE II: Comparison between three choices of the weighing factor.

| Optimum w.r.t | $\mathcal{F}_{\min}$ change w.r.t. case 1 (%) | $\mathcal{T}_{\min}$ change w.r.t. case 1 (%) | $\mathcal{F}_{\min}$ change w.r.t. case 2 (%) | $\mathcal{T}_{\min}$ change w.r.t. case 2 (%) |
|---|---|---|---|---|
| Case 1 | - | - | +41.42 | -18.35 |
| Case 2 | -29.29 | +22.47 | - | - |
| Case 3 | -4.38 | +12.58 | +35.23 | -8.08 |

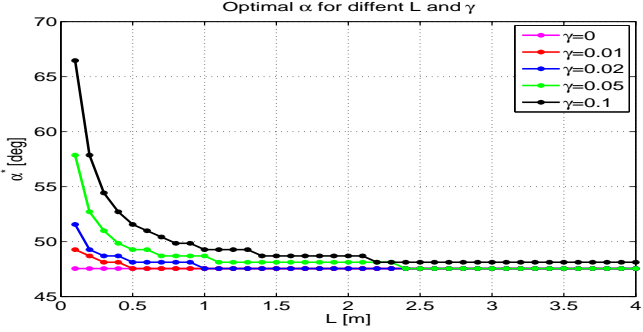

Fig. 6: Optimal design angles $\alpha^*$ for different UAV scales and propeller-specific drag-to-thrust ratios.

## III. DYNAMIC MODEL AND MATHEMATICAL PRELIMINARIES

In this section, the dynamic model of the fully actuated hexarotor will be presented. This is achieved by considering it as a rigid body in the special Euclidean group $SE(3) := SO(3) \ltimes \mathbb{R}^3$. In fact, the model presented is generic and represents any fully actuated aerial vehicle that can be modeled as a floating rigid body in space.

Let $\boldsymbol{\xi}_B^I \in \mathbb{R}^3$ represent the Cartesian position of the origin of the body fixed frame $o_B$ in $\Psi_I$, while $\boldsymbol{R}_B^I \in SO(3)$ represent the orientation of $\Psi_B$ with respect to $\Psi_I$. Let $\dot{\boldsymbol{\xi}}_B^I \in \mathbb{R}^3$ represent the linear velocity vector of the origin of $\Psi_B$ with respect to $\Psi_I$ expressed in $\Psi_I$, while $\boldsymbol{\omega}_B^{B,I} \in \mathbb{R}^3$ represents the angular velocity vector of $\Psi_B$ with respect to $\Psi_I$, expressed in $\Psi_B$. Let $m$ denote the mass of the vehicle, and $\boldsymbol{J} \in \mathbb{R}^{3\times 3}$ denote its constant mass moment of inertia matrix expressed in $\Psi_B$. The equations of motion of the vehicle are given by the *Euler-Poincaré* equations

$$\dot{\boldsymbol{H}}_B^I = \boldsymbol{H}_B^I \tilde{\boldsymbol{T}}_B^{B,I}, \tag{14}$$

$$\mathcal{I}\dot{\boldsymbol{T}}_B^{B,I} = \boldsymbol{ad}_{\boldsymbol{T}_B^{B,I}}^\top (\mathcal{I}\boldsymbol{T}_B^{B,I}) + \boldsymbol{W}_g^B + \boldsymbol{W}_p^B, \tag{15}$$

where $\boldsymbol{H}_B^I$ denotes a $4 \times 4$ homogenous matrix that represents coordinate transformation from $\Psi_B$ to $\Psi_I$, defined as

$$\boldsymbol{H}_B^I := \begin{bmatrix} \boldsymbol{R}_B^I & \boldsymbol{\xi}_B^I \\ \boldsymbol{0} & 1 \end{bmatrix} \in SE(3). \tag{16}$$

The twist $\tilde{\boldsymbol{T}}_B^{B,I}$ of $\Psi_B$ with respect to $\Psi_I$, expressed in $\Psi_B$ is defined as

$$\tilde{\boldsymbol{T}}_B^{B,I} := \begin{bmatrix} \tilde{\boldsymbol{\omega}}_B^{B,I} & \boldsymbol{v}_B^{B,I} \\ \boldsymbol{0} & 0 \end{bmatrix} := \begin{bmatrix} \boldsymbol{R}_I^B \dot{\boldsymbol{R}}_B^I & \boldsymbol{R}_I^B \dot{\boldsymbol{\xi}}_B^I \\ \boldsymbol{0} & 0 \end{bmatrix} \in se(3), \tag{17}$$

which is an element of the Lie algebra of $SE(3)$. The *tilde map* (operating on $\boldsymbol{\omega}$) denotes the skew-symmetric matrix, i.e., $(\cdot)^\sim : \mathbb{R}^3 \to so(3)$, defined such that $\tilde{\boldsymbol{\omega}}\boldsymbol{x} = \boldsymbol{\omega} \wedge \boldsymbol{x}, \forall \boldsymbol{x} \in \mathbb{R}^3$. The vector of Plücker coordinates of the twist is denoted by $\boldsymbol{T}_B^{B,I}$, and defined as

$$\boldsymbol{T}_B^{B,I} := \begin{bmatrix} \boldsymbol{\omega}_B^{B,I} \\ \boldsymbol{v}_B^{B,I} \end{bmatrix} \in \mathbb{R}^6. \tag{18}$$

The symbol $\mathcal{I}$ denotes the generalized inertia tensor of the UAV expressed in $\Psi_B$ and is represented by

$$\mathcal{I} := \begin{bmatrix} \boldsymbol{J} & \boldsymbol{0} \\ \boldsymbol{0} & m\boldsymbol{I}_3 \end{bmatrix} \in \mathbb{R}^{6\times 6}, \tag{19}$$

where $\boldsymbol{I}_i$ denotes the identity matrix of dimension $i$, and $\boldsymbol{0}$ denotes a matrix with all entries equal to zero. The mapping $\boldsymbol{ad}_{\boldsymbol{T}} : se(3) \to se(3)$ denotes the adjoint representation of the Lie algebra, which is a function of the twist $\boldsymbol{T}$, expressed as

$$\boldsymbol{ad}_{\boldsymbol{T}} := \begin{bmatrix} \tilde{\boldsymbol{\omega}} & \boldsymbol{0} \\ \tilde{\boldsymbol{v}} & \tilde{\boldsymbol{\omega}} \end{bmatrix}. \tag{20}$$

Moreover, the wrenches applied to the rigid body include the propellers' control wrench $\boldsymbol{W}_p^B \in \mathbb{R}^6$, previously defined in (5), and the gravity wrench $\boldsymbol{W}_g^B \in \mathbb{R}^6$ expressed as

$$\boldsymbol{W}_g^B = \boldsymbol{Ad}_{\boldsymbol{H}_B^{\mathring{B}}}^\top [0,0,0,0,0,-mg]^\top, \tag{21}$$

where the mapping $\boldsymbol{Ad}_{\boldsymbol{H}} : se(3) \to se(3)$ denotes the adjoint representation of the Lie group $SE(3)$, given by

$$\boldsymbol{Ad}_{\boldsymbol{H}} := \begin{bmatrix} \boldsymbol{R} & \boldsymbol{0} \\ \tilde{\boldsymbol{\xi}}\boldsymbol{R} & \boldsymbol{R} \end{bmatrix}, \tag{22}$$

and the homogenous matrix $\boldsymbol{H}_B^{\mathring{B}}$ represents the coordinate transformation from the body-fixed frame $\Psi_B$ to another body-fixed frame $\Psi_{\mathring{B}}$ which is attached to the CoM but aligned with the inertial reference frame, i.e.,

$$\boldsymbol{H}_B^{\mathring{B}} = \begin{bmatrix} \boldsymbol{R}_B^I & \boldsymbol{0} \\ \boldsymbol{0} & 1 \end{bmatrix}. \tag{23}$$

In addition, it is assumed that gravity is downwards parallel to $\hat{\boldsymbol{z}}_I$, with $g$ denoting the gravitational acceleration constant.

The following identities are useful in the control system design in section IV. Let $\Psi_j, \Psi_k, \Psi_m$ denote three different reference frames. Then, the following identities hold

$$\boldsymbol{T}_j^{j,k} = -\boldsymbol{T}_k^{j,j} \tag{24}$$

$$\tilde{\boldsymbol{T}}_j^{m,k} = \boldsymbol{H}_j^m \tilde{\boldsymbol{T}}_j^{j,k} \boldsymbol{H}_m^j \tag{25}$$

$$\boldsymbol{T}_j^{m,k} = \boldsymbol{Ad}_{\boldsymbol{H}_j^m} \boldsymbol{T}_j^{j,k} \tag{26}$$

$$\frac{d}{dt}(\boldsymbol{Ad}_{\boldsymbol{H}_j^k}) = \boldsymbol{Ad}_{\boldsymbol{H}_j^k} \boldsymbol{ad}_{\boldsymbol{T}_j^{j,k}} \tag{27}$$

$$(\boldsymbol{Ad}_{\boldsymbol{H}_j^k})^{-1} = \boldsymbol{Ad}_{\boldsymbol{H}_k^j} \tag{28}$$

$$\boldsymbol{ad}_{\boldsymbol{T}_j^{m,k}} = \boldsymbol{Ad}_{\boldsymbol{H}_j^m} \boldsymbol{ad}_{\boldsymbol{T}_j^{j,k}} \boldsymbol{Ad}_{\boldsymbol{H}_m^j} \tag{29}$$

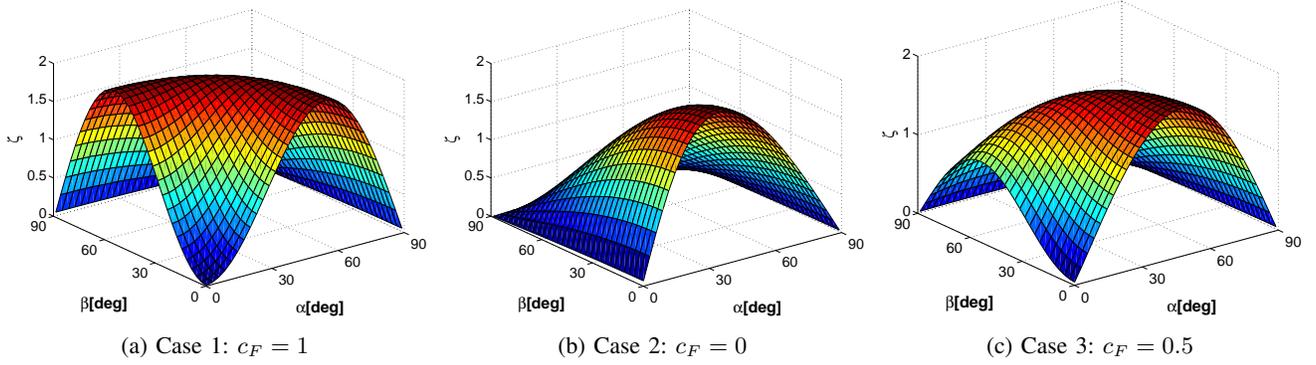

(a) Case 1: $c_F = 1$     (b) Case 2: $c_F = 0$     (c) Case 3: $c_F = 0.5$

Fig. 5: Visualization of three choices for the objective function defined in (13). The vertical axis displays the dimensionless value $\zeta$ of the objective function.

Proofs can be found in [17].

## IV. Geometric Tracking Controller Design

In this section, we present the geometric control system design for the class of fully actuated aerial vehicles, represented by (14,15). The control problem for the class of UAVs under study is to design a control law for the propellers' thrusts $\boldsymbol{\lambda}$ that enables the tracking of arbitrary trajectories for the UAV's configuration $\boldsymbol{H}_B^I$. In order for the controller design methodology to be generic, it will be assumed that the desired control wrench $\boldsymbol{W}_{p,\text{des}}^B$ on the vehicle's body results from the control law. Then, in a model-based manner, the desired propellers' thrust are computed by the inverse of the mapping $\boldsymbol{M}$ given by

$$\boldsymbol{\lambda}_{\text{des}} = \boldsymbol{M}^{-1} \boldsymbol{W}_{p,\text{des}}^B, \tag{30}$$

where the map $\boldsymbol{M}$ is platform specific.

The methodology used in this work for the control system design is based on the framework introduced by Bullo and Murray [14] for geometric tracking on SE(3). Let $\boldsymbol{H}_D^I \in SE(3)$ denote the desired frame to be tracked by the UAV, and $\boldsymbol{T}_D^{D,I}$ denote the desired twist expressed in the desired reference frame $\Psi_D$, such that $\tilde{\boldsymbol{T}}_D^{D,I} := \boldsymbol{H}_I^D \dot{\boldsymbol{H}}_D^I$. The desired trajectory is assumed to be at least twice differentiable. For the readability of the remaining sections, the following notations are used henceforth for any frame $\Psi_j$; $\boldsymbol{H}_j := \boldsymbol{H}_j^I$, and $\boldsymbol{T}_j := \boldsymbol{T}_j^{j,I}$, and the same applies for their components.

Let $\boldsymbol{H}_{e,r} \in SE(3)$ denote the right group error which describes the relative configuration from $\Psi_B$ to $\Psi_D$, expressed as

$$\boldsymbol{H}_{e,r} := \boldsymbol{H}_B^D = \boldsymbol{H}_D^{-1} \boldsymbol{H}_B = \begin{bmatrix} \boldsymbol{R}_{e,r} & \boldsymbol{\xi}_e \\ \boldsymbol{0} & 1 \end{bmatrix}, \tag{31}$$

where $\boldsymbol{R}_{e,r}$ represents the right orientation tracking error on $SO(3)$, and $\boldsymbol{\xi}_e$ represents the position tracking error, expressed in the desired frame. The terms $\boldsymbol{R}_{e,r}$ and $\boldsymbol{\xi}_e$ are given by

$$\boldsymbol{R}_{e,r} := \boldsymbol{R}_D^\top \boldsymbol{R}_B \tag{32}$$

$$\boldsymbol{\xi}_e := \boldsymbol{R}_D^\top \bar{\boldsymbol{\xi}}_e := \boldsymbol{R}_D^\top (\boldsymbol{\xi}_B - \boldsymbol{\xi}_D). \tag{33}$$

A compatible twist error with the configuration error $\boldsymbol{H}_{e,r}$ is the right twist error $\tilde{\boldsymbol{T}}_{e,r}$ defined by

$$\tilde{\boldsymbol{T}}_{e,r} := \boldsymbol{H}_{e,r}^{-1} \dot{\boldsymbol{H}}_{e,r}, \tag{34}$$

which in vector form can be expressed as

$$\boldsymbol{T}_{e,r} = \boldsymbol{T}_B - Ad_{\boldsymbol{H}_{e,r}^{-1}} \boldsymbol{T}_D. \tag{35}$$

The detailed proof of (35) can be found in the appendix.

To proceed with the controller design, a measure of the configuration error $\boldsymbol{H}_{e,r}$ is essential. In this work, we utilize the notion of an error function, which denotes a positive definite map on the Lie group $\phi : SE(3) \to \mathbb{R}$, introduced in [14]. By composing different group errors on $SE(3)$ with a chosen positive definite function, there are many possibilities for the error function [14]. The error function utilized in this work is given by

$$\phi(\boldsymbol{H}_{e,r}) := \frac{1}{2} \operatorname{tr}(\boldsymbol{K}_{p_1}(\boldsymbol{I}_3 - \boldsymbol{R}_{e,r})) + \frac{1}{2} k_{p_2} \|\boldsymbol{\xi}_e\|^2, \tag{36}$$

where $\boldsymbol{K}_{p_1} > 0 \in \mathbb{R}^{3\times 3}$ and $k_{p_2} > 0 \in \mathbb{R}$ are the controller's proportional gains, and $\|\cdot\|$ denotes the 2-norm. This choice for the error function yields a coordinate-free controller that is invariant under a change of inertial frame [14].

The time derivative of the error function (36) evaluated on the configuration error trajectory, is given by

$$\dot{\phi}(\boldsymbol{H}_{e,r}, \dot{\boldsymbol{H}}_{e,r}) := -\frac{1}{2} \operatorname{tr}(\boldsymbol{K}_{p_1} \dot{\boldsymbol{R}}_{e,r}) + k_{p_2} \boldsymbol{\xi}_e^\top \dot{\boldsymbol{\xi}}_e. \tag{37}$$

With some mathematical manipulation, we can factor out the twist error $\boldsymbol{T}_{e,r}$, and then express (37) as

$$\dot{\phi} = [(\operatorname{as}(\boldsymbol{K}_{p_1} \boldsymbol{R}_{e,r})^\vee)^\top, \ -k_{p_2} \boldsymbol{\xi}_e^\top \boldsymbol{R}_{e,r}] \boldsymbol{T}_{e,r}, \tag{38}$$

where $\operatorname{as}(\cdot)$ denotes the skew-symmetric part of a given matrix, i.e., $\operatorname{as}(\boldsymbol{A}) = \frac{1}{2}(\boldsymbol{A} - \boldsymbol{A}^\top)$, and the *vee map* $(\cdot)^\vee : so(3) \to \mathbb{R}^3$ is the inverse of the *tilde* map. The derivation of (38) can be found in the appendix.

Now we design the actual control law of (15) by deriving the dynamics of the error twist $\boldsymbol{T}_{e,r}$, given by

$$\dot{\boldsymbol{T}}_{e,r} = \dot{\boldsymbol{T}}_B - Ad_{\boldsymbol{H}_{e,r}^{-1}} \dot{\boldsymbol{T}}_D - (\dot{Ad}_{\boldsymbol{H}_{e,r}^{-1}}) \boldsymbol{T}_D. \tag{39}$$

By substituting (15) in (39), the twist error dynamics can be written as

$$\dot{T}_{e,r} = \mathcal{I}^{-1}(ad_{T_B}^\top(\mathcal{I}T_B) + W_g^B + W_p^B) \\ - Ad_{H_{e,r}^{-1}}\dot{T}_D - (\dot{Ad}_{H_{e,r}^{-1}})T_D. \quad (40)$$

Finally, we design the control wrench $W_p^B$ to be

$$W_p^B = W_{ff} + W_p + W_d, \quad (41)$$

where $W_{ff}$ is a feedforward term and $W_p, W_d$ are the proportional and derivative feedback terms, respectively. The feedforward term is given by

$$W_{ff} = -W_g^B - ad_{T_B}^\top \mathcal{I}T_B \\ + \mathcal{I}(Ad_{H_{e,r}^{-1}}\dot{T}_D - ad_{T_{e,r}}Ad_{H_{e,r}^{-1}}T_D). \quad (42)$$

The proportional feedback term is given by

$$W_p = -\begin{bmatrix} \text{as}(K_{p_1}R_{e,r})^\vee \\ -k_{p_2}R_{e,r}^\top \xi_e \end{bmatrix}, \quad (43)$$

which is extracted from (38), while the derivative feedback term is given by

$$W_d = -K_d T_{e,r}, \quad (44)$$

with $K_d \in \mathbb{R}^{6\times 6}$ denoting a positive definite matrix of controller gains. The last term in the feedforwad term (42) is equivalent to the last term in the right hand side of (40); the proof can be found in the appendix. The stability properties of the closed loop system can be found in [14, Lemma 10] or in [19, Theorem 11.29].

## V. SIMULATION RESULTS

In this section, the simulation results of the geometric tracking controller (41-44) implemented on the omnidirectional hexarotor will be presented. The hexarotor simulation parameters are as follows: $m = 0.6$ kg, $g = 9.81$ m/s$^2$, and $J = \text{diag}(0.5, 0.5, 2)$. The controller gains are chosen to be $K_{p_1} = 10I_3$, $k_{p_2} = 3$, and $K_d = \text{diag}(5,5,5,2,2,2)$.

The UAV is simulated to track a desired position and orientation trajectory from a non-hovering initial configuration. The desired position trajectory is defined by $\xi_D = [r_d\cos(w_d t), r_d\sin(w_d t), r_d\sin(w_d t) + z_{off}]^\top$ where $r_d, w_d, z_{off}$ characterize the trajectory's shape. The desired orientation trajectory is defined by $\dot{\omega}_D = [\sin(t), 2\sin(t), \sin(t)]^\top$ with $\omega_D$ and $R_D$ computed by integration. The desired configuration, desired twist, and its derivative, respectively, are computed by

$$H_D = \begin{bmatrix} R_D & \xi_D \\ 0 & 1 \end{bmatrix}, T_D = \begin{bmatrix} \omega_D \\ R_D^\top \dot{\xi}_D \end{bmatrix},$$

$$\dot{T}_D = \begin{bmatrix} \dot{\omega}_D \\ R_D^\top \ddot{\xi}_D - \tilde{\omega}_D R_D^\top \dot{\xi}_D \end{bmatrix}.$$

Moreover, at $t = 7$ s, an impulse disturbance torque is applied in simulation to validate the robustness of the controller.

The simulation results can be found in Figs. 7-9. In Fig. 7, the configuration error function (36) is displayed. For the comprehensibility of the results, the rotational and translational components of the error function have been plotted separately. Moreover, in Fig. 8, the elements of the error twist (35) are presented. The first three elements represent the rotational part of the twist, which successfully recovers from the disturbance applied at $t = 7$ s. The twist errors are initially non-zero due to the non-hovering initial state of the UAV. The convergence of both the configuration and twist tracking errors is validated through this simulation. Moreover, the control effort required to achieve this tracking performance is presented in Fig 9. A video showing the simulation results and the tracking performance can be found in the supplementary material.

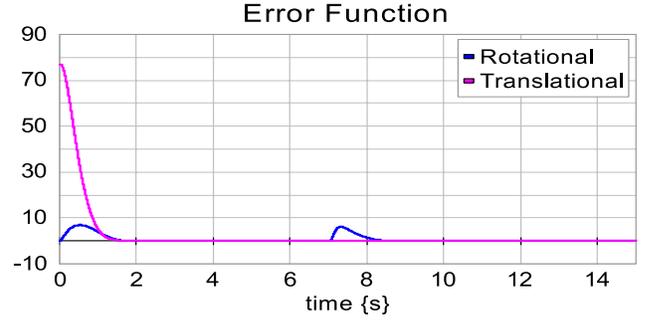

Fig. 7: Translational and rotational components of the error function which measures the size of the configuration error between the desired and actual reference frames. At $t = 7$ s, an impulse torque disturbance is applied.

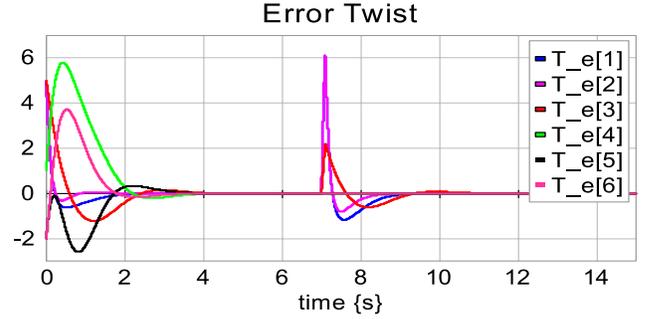

Fig. 8: Elements of the error twist compatible with the chosen configuration error. The first three elements comprise the rotational part of the twist, while the other three comprise the translational part.

## VI. CONCLUSION AND FUTURE WORK

In this work, we presented the design, modeling, and control of an omnidirectional hexarotor, optimized for interaction. Based on the criteria of maximizing the control wrench, it was shown that only axial tilting of the propellers with respect to the UAV's body yields optimal results. It was shown that the optimal design angle is invariant to the thrust-to-drag ratio of the propellers for large-scale UAVs. In contrast to the conventional hexarotor, the proposed hexarotor can generate at least 1.9 times the maximum thrust of one rotor in any direction, in addition to the higher control torque around the vehicle's upward axis. A geometric tracking controller on $SE(3)$, that is invariant to inertial coordinate changes, was proposed for the class of fully actuated UAVs and validated in simulation. Future work includes the experimental validation of the proposed UAV design and control system. Moreover, we intend to extend the proposed geometric controller for physical interaction with the environment.

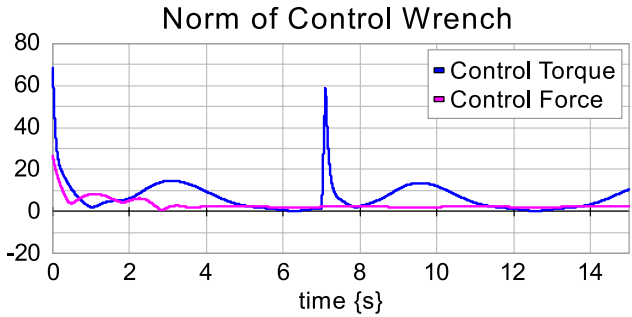

Fig. 9: The control wrench applied to the UAV's body to track the specified trajectory. The norm of the control torques is plotted separately from the norm of the control forces.

## APPENDIX

In this appendix, the proofs of (35), (38), and the last term in (42) will be presented, respectively. Starting from the definition of the right twist error given in (34), and by substituting in it (31), we get

$$\tilde{T}_{e,r} = H_{e,r}^{-1}\dot{H}_{e,r},$$
$$= (H_I^D H_B^I)^{-1}\frac{d}{dt}(H_I^D H_B^I),$$
$$= H_I^B H_D^I H_I^D \dot{H}_B^I + H_I^B H_D^I \dot{H}_I^D H_B^I,$$
$$= H_I^B \dot{H}_B^I + H_I^B \tilde{T}_I^{I,D} H_B^I,$$
$$= \tilde{T}_B^{B,I} - H_I^B \tilde{T}_D^{I,I} H_B^I,$$
$$= \tilde{T}_B^{B,I} - \tilde{T}_D^{B,I},$$

where in the last two steps the identities (24) and (25) were used. Furthermore, by using identity (26), we can write the vector form of $\tilde{T}_{e,r}$ as $T_{e,r} = T_B^{B,I} - Ad_{H_D^B} T_D^{D,I}$. Then, by the definition of $H_{e,r}$ in (31), we get (35).

For the proof of the last term in (42), we start by using (31) and identity (27) such that

$$(\dot{Ad}_{H_{e,r}^{-1}}) = (\dot{Ad}_{H_D^B}) = Ad_{H_D^B} ad_{T_D^{D,B}},$$
$$= Ad_{H_D^B} Ad_{H_B^D} ad_{T_D^{B,B}} Ad_{H_D^B},$$
$$= ad_{T_D^{B,B}} Ad_{H_D^B},$$
$$= -ad_{T_B^{B,D}} Ad_{H_D^B},$$

where identities (29), (28), and (24) where used respectively. Finally, by the definition of $T_{e,r} := T_B^{B,D}$ and (31), we get the proposed form in (42).

For the proof of (38), it consists of two parts; the first is proving that

$$-\frac{1}{2}\mathrm{tr}(K_{p_1}\dot{R}_{e,r}) \equiv (\mathrm{as}(K_{p_1}R_{e,r})^{\vee})^{\top}\omega_{e,r}, \quad (45)$$

while the second part is to prove that

$$\xi_e^{\top}\dot{\xi}_e \equiv \xi_e^{\top} R_{e,r} v_{e,r}. \quad (46)$$

The terms $\omega_{e,r}$ and $v_{e,r}$ are the components of the twist $T_{e,r}$. By expanding Eq. (35) using (22) and (31), we get

$$T_{e,r} = \begin{bmatrix} \omega_{e,r} \\ v_{e,r} \end{bmatrix} = \begin{bmatrix} \omega_B - R_{e,r}^{\top}\omega_D \\ v_B - R_{e,r}^{\top}(v_D + \tilde{\omega}_D \xi_e) \end{bmatrix},$$

where $\omega_k, v_k$ are the corresponding elements of the twist $T_k$.

Now we start with the first term (45), which is proven as follows:

$$-\frac{1}{2}\mathrm{tr}(K_{p_1}\dot{R}_{e,r}),$$
$$= -\frac{1}{2}\mathrm{tr}(K_{p_1}R_{e,r}R_{e,r}^{\top}\dot{R}_{e,r}),$$
$$= -\frac{1}{2}\mathrm{tr}(K_{p_1}R_{e,r}\tilde{\omega}_{e,r}),$$
$$= -\frac{1}{2}\mathrm{tr}([\mathrm{as}(K_{p_1}R_{e,r}) + \mathrm{sym}(K_{p_1}R_{e,r})]\tilde{\omega}_{e,r}),$$
$$= -\frac{1}{2}\mathrm{tr}(\mathrm{as}(K_{p_1}R_{e,r})\tilde{\omega}_{e,r}),$$
$$= (\mathrm{as}(K_{p_1}R_{e,r})^{\vee})^{\top}\omega_{e,r},$$

where we used the identity $\mathrm{tr}(\tilde{x}\tilde{y}) = -2x^{\top}y$, and the orthogonality property of the subspaces of symmetric and antisymmetric matrices with respect to the trace inner product [19].

For the second part (46), we start by substituting (33) in the left hand side of (46), such that

$$\xi_e^{\top}\dot{\xi}_e = (R_D^{\top}\bar{\xi}_e)^{\top}\frac{d}{dt}(R_D^{\top}\bar{\xi}_e),$$
$$= \bar{\xi}_e^{\top} R_D R_D^{\top}\dot{\bar{\xi}}_e + \bar{\xi}_e^{\top} R_D \dot{R}_D^{\top}\bar{\xi}_e,$$
$$= \bar{\xi}_e^{\top} R_D R_D^{\top}(\dot{\xi}_B - \dot{\xi}_D) + \bar{\xi}_e^{\top}\tilde{\omega}_I^{I,D}\bar{\xi}_e,$$
$$= \bar{\xi}_e^{\top}(R_D R_D^{\top}\dot{\xi}_B - R_D R_D^{\top}\dot{\xi}_D) - \bar{\xi}_e^{\top}\tilde{\omega}_D^{D,I}\bar{\xi}_e,$$
$$= \bar{\xi}_e^{\top}(\dot{\xi}_B - R_D R_D^{\top}\dot{\xi}_D) - \bar{\xi}_e^{\top} R_D \tilde{\omega}_D^{D,I} R_D^{\top}\bar{\xi}_e,$$
$$= \bar{\xi}_e^{\top}(\dot{\xi}_B - R_D v_D) - \bar{\xi}_e^{\top} R_D \tilde{\omega}_D^{D,I}\xi_e,$$
$$= \bar{\xi}_e^{\top} R_B R_B^{\top}(\dot{\xi}_B - R_D v_D - R_D \tilde{\omega}_D^{D,I}\xi_e),$$
$$= \bar{\xi}_e^{\top} R_B(R_B^{\top}\dot{\xi}_B - R_B^{\top} R_D v_D - R_B^{\top} R_D \tilde{\omega}_D^{D,I}\xi_e),$$
$$= \bar{\xi}_e^{\top} R_B(v_B - R_{e,r}^{\top} v_D - R_{e,r}^{\top}\tilde{\omega}_D^{D,I}\xi_e),$$
$$= \bar{\xi}_e^{\top} R_B v_{e,r},$$
$$= \xi_e^{\top} R_{e,r} v_{e,r},$$

which concludes the proof.

## ACKNOWLEDGMENT

The authors would like to thank Dannis Brouwer and Ronald Aarts for their contribution in the optimization-based design.

## REFERENCES

[1] H. Yang and D. Lee, "Dynamics and control of quadrotor with robotic manipulator," in *Robotics and Automation (ICRA), 2014 IEEE International Conference on*. IEEE, 2014, pp. 5544–5549.

[2] T. Bartelds, A. Capra, S. Hamaza, S. Stramigioli, and M. Fumagalli, "Compliant aerial manipulators: Toward a new generation of aerial robotic workers," *IEEE Robotics and Automation Letters*, vol. 1, no. 1, pp. 477–483, 2016.

[3] S. Salazar, H. Romero, R. Lozano, and P. Castillo, "Modeling and real-time stabilization of an aircraft having eight rotors," in *Unmanned Aircraft Systems*. Springer, 2008, pp. 455–470.

[4] M. Ryll, H. H. Bulthoff, and P. R. Giordano, "Modeling and control of a quadrotor UAV with tilting propellers," in *2012 IEEE Int. Conf. Robot. Autom.* IEEE, may 2012, pp. 4606–4613.

[5] M. Ryll, H. H. Bülthoff, and P. R. Giordano, "A novel overactuated quadrotor unmanned aerial vehicle: Modeling, control, and experimental validation," *IEEE Trans. Control Syst. Technol.*, vol. 23, no. 2, pp. 540–556, 2015.

[6] R. Voyles and G. Jiang, "Hexrotor UAV platform enabling dextrous interaction with structures-Preliminary work," in *Safety, Security, and Rescue Robotics (SSRR), 2012 IEEE International Symposium on*. IEEE, 2012, pp. 1–7.

[7] K. Kiso, T. Ibuki, M. Yasuda, and M. Sampei, "Structural optimization of hexrotors based on dynamic manipulability and the maximum translational acceleration," in *2015 IEEE Conf. Control Appl.* IEEE, sep 2015, pp. 774–779.

[8] S. Rajappa, M. Ryll, H. H. Bulthoff, and A. Franchi, "Modeling, control and design optimization for a fully-actuated hexarotor aerial vehicle with tilted propellers," *Proc. - IEEE Int. Conf. Robot. Autom.*, vol. 2015-June, no. June, pp. 4006–4013, 2015.

[9] M. Ryll, G. Muscio, F. Pierri, E. Cataldi, G. Antonelli, F. Caccavale, and A. Franchi, "6D physical interaction with a fully actuated aerial robot," in *2017 IEEE International Conference on Robotics and Automation*, 2017.

[10] D. Brescianini and R. D'Andrea, "Design, modeling and control of an omni-directional aerial vehicle," in *2016 IEEE Int. Conf. Robot. Autom.* IEEE, may 2016, pp. 3261–3266.

[11] S. Park, J. Her, J. Kim, and D. Lee, "Design, modeling and control of omni-directional aerial robot," in *2016 IEEE/RSJ Int. Conf. Intell. Robot. Syst.* IEEE, oct 2016, pp. 1570–1575.

[12] T. Lee, M. Leoky, and N. H. McClamroch, "Geometric tracking control of a quadrotor UAV on SE(3)," in *Decision and Control (CDC), 2010 49th IEEE Conference on*. IEEE, 2010, pp. 5420–5425.

[13] R. Mahony, V. Kumar, and P. Corke, "Multirotor aerial vehicles," *IEEE Robotics and Automation magazine*, vol. 20, no. 32, pp. 20–32, 2012.

[14] F. Bullo and R. M. Murray, "Tracking for fully actuated mechanical systems: a geometric framework," *Automatica*, vol. 35, no. 1, pp. 17–34, 1999.

[15] M. Cutler, N. K. Ure, B. Michini, and J. P. How, "Comparison of fixed and variable pitch actuators for agile quadrotors," in *AIAA Guidance, Navigation, and Control Conference (GNC)*, vol. 2, 2011.

[16] P. Bosscher, A. T. Riechel, and I. Ebert-Uphoff, "Wrench-feasible workspace generation for cable-driven robots," *IEEE Transactions on Robotics*, vol. 22, no. 5, pp. 890–902, 2006.

[17] S. Stramigioli, *Modeling and IPC control of interactive mechanical system - A coordinate-free approach*. Springer-Verlag London, 2001.

[18] R. Rashad, P. Kuipers, J. Engelen, and S. Stramigioli, "Design, Modeling and Geometric Control on SE(3) of a Fully-Actuated Hexarotor for Aerial Interaction," in *arXiv*, 2017.

[19] F. Bullo and A. D. Lewis, *Geometric control of mechanical systems: modeling, analysis, and design for simple mechanical control systems*. Springer Science & Business Media, 2004, vol. 49.